\documentclass[runningheads,a4paper]{llncs}
\usepackage{latexsym,amsfonts,amssymb}
\usepackage{graphicx}
\usepackage{mathtools}

\begin{document}
	\mainmatter  
	
	\title{Nilpotent approximation of a trident snake robot controlling distribution}
	
	\titlerunning{Nilpotent approximation of a trident snake robot controlling distribution}
	
	%
	%
	\author{Jaroslav Hrdina, Ale\v s N\'avrat and Petr Va\v{s}\'{i}k
		\thanks{The research
			was supported by a grant no. FSI-S-14-2290.} }
	\authorrunning{J. Hrdina, A. N\'avrat, P. Va\v{s}\'{i}k}
	
	\institute{
		Brno University of Technology, Faculty of Mechanical Engineering\\
		Institute of Mathematics\\
		Technick\'{a} 2896/2, 616 69 Brno\\
		Czech Republic}
	
	%
	%
	
	\maketitle

	\begin{abstract}
		We construct a privileged system of coordinates with respect to the controlling distribution of a trident snake robot and, furthermore, we construct a nilpotent approximation with respect to the given filtration. Note that all constructions are local in the neighbourhood of a particular point. We compare the motions corresponding to the Lie bracket of the original controlling vector fields and their nilpotent approximation.
	\end{abstract}
	
	\section{Introduction}
	Within this paper, we consider a trident snake robot moving on a 
	planar surface. More precisely, it is a model when to each vertex of an equilateral triangle a leg of length 1 is attached that is endowed by a pair 
	of passive wheels at its end. The joints of the legs to the triangle platform are motorised and thus the possible motion directions are 
	determined uniquely. Local controllability of such mechanism is known, see \cite{Ishikawa}.
	If the generalized coordinates are considered, the non--holonomic 
	forward kinematic equations can be understood as a Pfaff system and its 
	solution as a distribution in the configuration space. Rachevsky--Chow 
	Theorem implies that the appropriate non--holonomic system is locally 
	controllable if the corresponding distribution is not integrable and the 
	span of the Lie algebra generated by the controlling distribution has to 
	be of the same dimension as the configuration space. 
	The spanned Lie algebra is then naturally endowed by a filtration which 
	shows the way to realize the motions by means of the vector field 
	brackets \cite{b4,b5}. In our case, the system is locally controllable and the filtration is $(3,6)$.
	
	In order to simplify the trident snake robot control, in Section \ref{Privcoord} we construct a privileged system of coordinates with respect to the distribution given by local nonholonomic conditions and, furthermore, in Section \ref{NilApprox} we construct a nilpotent approximation of the transformed distribution with respect to the given filtration.. Note that all constructions are local in the neighbourhood of 0.
	
	Finally, we compare the motions generated by the Lie brackets of the original controlling vector fields and their nilpotent approximation. The accuracy is demonstrated by simulations in MATLAB.

	\section{Preliminaries}\label{Prelim}
	We recall the following concepts of functions or vector fields orders and distribution weights, see \cite{Jean}.
	Let $X_1,...,X_m$ denote the smooth vector fields on a manifold $M$ and $C^\infty (p)$ denote the set of germs of smooth functions at $p\in M$. For $f\in C^\infty (p)$ we say that the Lie derivatives $X_if,X_iX_jf,...$ are non--holonomic derivatives of $f$ of order 1,2,... The non--holonomic derivative of order 0 of $f$ at $p$ is $f(p).$ 
	
	\begin{definition}
	Let $f\in C^\infty (p).$ Then the non--holonomic order of $f$ at $p$, denoted by $\mathrm{ord}_p(f)$, is the biggest integer $k$ such that all non--holonomic derivatives of $f$ of order smaller than $k$ vanish at $p.$
	\end{definition}
	Note that in case $M=\mathbb R^n,\ m=n$ and $X_i=\partial_{x_i},$ for a smooth function $f$, $\mathrm{ord}_0(f)$ is the smallest degree of monomials having nonzero coefficient in the Taylor series. In the language of non--holonomic derivatives, the order of a smooth function is given by the formula, \cite{Jean}:
	$$\mathrm{ord}_p(f)=\min \biggl\{ s\in\mathbb N: \exists i_1,...,i_s\in\{1,...,m\} \text{ s.t. } (X_{i_1}\cdots X_{i_s}f)(p)\neq 0 \biggr\} ,$$
	where the convention reads that $\min \emptyset = \infty.$ 
	
	If we denote by $\mathrm{VF}(p)$ the set of germs of smooth vector fields at $p\in M$, the notion of non--holonomic order extends to the vector fields as follows:
	
	\begin{definition}
		Let $X\in\mathrm{VF}(p).$ The non--holonomic order of $X$ at $p$, denoted by $\mathrm{ord}_p(X),$ is a real number defined by:
		$$\mathrm{ord}_p(X)=\sup \biggl\{ \sigma\in\mathbb R : \mathrm{ord}_p(Xf)\geq\sigma+\mathrm{ord}_p(f), \forall f\in C^\infty (p)\biggr\}.$$
	\end{definition}
	
	Note that $\mathrm{ord}_p(X)\in\mathbb Z.$ Moreover, the null vector field $X\equiv 0$ has infinite order, $\mathrm{ord}_p(0)=\infty.$ Furthermore, $X_1,...,X_m$ are of order $\geq -1$, $[X_i,X_j]$ of order $\geq -2$, etc. 
	
	Using the notion of a vector field order one can define 
	
	\begin{definition}
		A family of $m$ vector fields $(\hat X_1,...,\hat X_m)$ defined near $p$ is called a first order approximation of $(X_1,...,X_m)$ at $p$ if the vector fields $X_i-\hat X_i, i=1,...,m$ are of order $\geq 0$ at $p.$
	\end{definition}
	
	Finally, to define the weights of distributions we use the same notation as in \cite{Jean}. Let us by $\varDelta^1$ denote the distribution 
	$$\varDelta^1=\mathrm{span}\{X_1,...,X_m\}$$
	and for $s\geq 1$ define 
	$$\varDelta^{s+1}=\varDelta^s+[\varDelta^1,\varDelta^s],$$
	where $[\varDelta^1,\varDelta^s]=\mathrm{span}\{[X,Y]:X\in\varDelta^1, Y\in\varDelta^s\}.$ Then
	$$\varDelta^s=\mathrm{span}\{X_I:|I|\leq s\}.$$
	Note that this directly leads to the fact that every $X\in\varDelta^s$ is of order $\geq -s.$ Now let us consider the sequence
	$$\varDelta^1(p)\subset\varDelta^2(p)\subset\cdots\subset\varDelta^{r-1}\varsubsetneq\varDelta^r(p)=T_pM,$$ where  $r=r(p)$ is called the degree of non--holonomy at $p.$ Set $n_i(p)=\dim\varDelta^i(p).$ Then we can define the weights at $p$, $w_i=w_i(p), i=1,...,n=n_{r(p)}$ by setting $w_j=s$ if $n_{s-1}(p)<j\leq n_s(p),$ where $n_0=0.$ In other words, we have 
	$$w_1=\cdots=w_{n_1}=1,\ w_{n_1+1}=\cdots=w_{n_2}=2,...,w_{n_{r-1}+1}=\cdots=w_{n_r}=r.$$
	The weights at $p$ form an increasing sequence $w_1(p)\leq\cdots\leq w_n(p).$
	
	\section{Trident snake robot}
	The mechanism of the trident snake robot was described in \cite{Ishikawa}. It consists of a body in the shape of an equilateral triangle with circumscribed circle of radius $r$ and three rigid links (also called legs) of constant length $l$ connected to the vertices of the triangular body by three motorised  joints. In this paper, we consider $r=1$ and $l=1$ . To each free link end, a pair of passive wheels is attached to provide an important snake-like property that the ground friction in the direction perpendicular to the link is considerably higher than the friction of a simple forward move.
	In particular, this prevents slipping sideways. 
	\begin{figure}[h]
		\centering
		\includegraphics[scale=0.5]{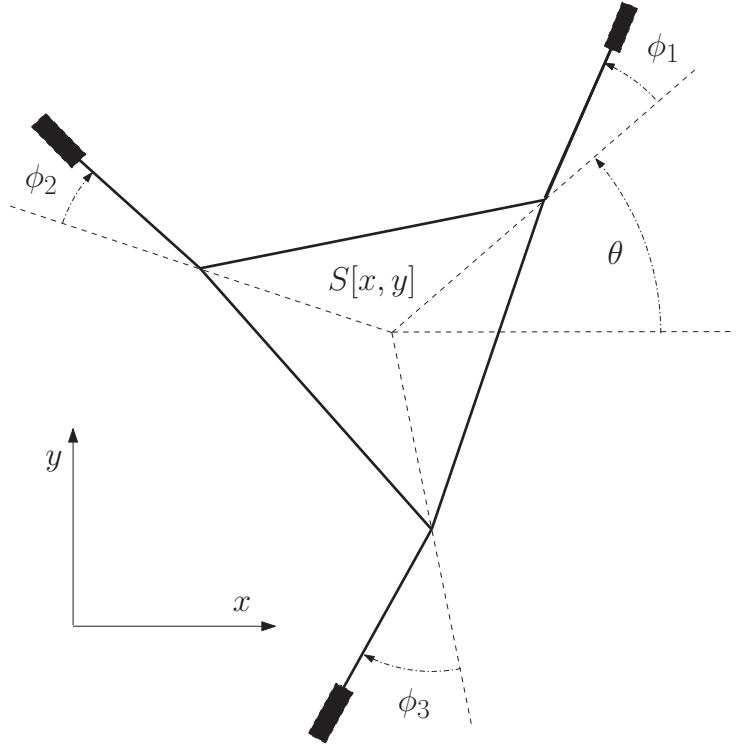}
		\caption{Trident snake robot model}
		\label{fig1}
	\end{figure}
	To describe the actual position of a trident snake robot we need the set of 6 generalized coordinates 
	$$q=(x,y,\theta,\phi_1,\phi_2.\phi_3)=:(x_1,x_2,x_3,x_4,x_5,x_6)$$ as shown in Figure 1.
	Hence the configuration space is (a subspace of) $\mathbb{R}^2\times S^1\times (S^1)^3$. Note that a fixed coordinate system $(x,y)$ is attached.

	\section{Local controllability and coordinate systems}
	Local controllability of such robot is given by the appropriate Pfaff system of ODEs. The solution with respect to $\dot q$ gives a control system $\dot{q}=G\mu$, where the control matrix $G$ is a $6\times 3$ matrix spanned by vector fields $g_1, g_2, g_3$, where
	\begin{align*}
	g_1 &=\cos\theta \partial_x +\sin\theta\partial_y+\sin\phi_1\partial_{\phi_1}+\sin(\phi_2+\tfrac{2\pi}{3})\partial_{\phi_2}
	+\sin(\phi_3+\tfrac{4\pi}{3})\partial_{\phi_3},
	\\
	g_2 &=-\sin\theta \partial_x +\cos\theta\partial_y-\cos\phi_1\partial_{\phi_1}-\cos(\phi_2+\tfrac{2\pi}{3})\partial_{\phi_2}
	-\cos(\phi_3+\tfrac{4\pi}{3})\partial_{\phi_3},
	\\
	g_3 &=\partial_\theta-(1+\cos\phi_1)\partial_{\phi_1}-(1+\cos\phi_2)\partial_{\phi_2}-(1+\cos\phi_3)\partial_{\phi_3}.
	\end{align*}
	Note that the parametrizations can vary by setting the angles within the triangular platform either $\tfrac{2\pi}{3}$ and $\tfrac{4\pi}{3}$ or $\tfrac{2\pi}{3}$ and $-\tfrac{2\pi}{3},$ etc.
	It is easy to check that in regular points these vector fields define a (bracket generating) distribution with growth vector $(3,6).$ It means that in each regular point the vector fields $g_1, g_2, g_3$ together with their Lie brackets span the whole tangent space.
	Consequently, the system is controllable by Chow--Rashevsky theorem.
	
	Let us decompose the control system in such way that the spatial coordinates $w:=(x,y,\theta)=(x_1,x_2,x_3)$ are parametrised by the angles $\phi:=(\phi_1,\phi_2.\phi_3)=(x_4,x_5,x_6)$, and, furthermore, the form where the invariant parameter $\theta$ is excluded, i.e. it is of the form
	\begin{equation}\label{eq1}A(\phi)R^T_\theta\dot w=\dot \phi,\end{equation}
	where 
	$$R^T_\theta=\begin{pmatrix}
	\cos\theta & \sin\theta & 0\\
	-\sin\theta & \cos\theta & 0\\
	0 & 0 & 1 
	\end{pmatrix}$$
	 is the matrix of rotation by the angle $\theta,$ see \cite{IshikawaMinami}. If the spatial coordinate transformation 
	\begin{equation}\label{eq2}v=(A(\phi))^{-1}\dot\phi\end{equation}
	is considered, we modify the system \eqref{eq1} and obtain
	$$\dot w=R_\theta (A(\phi))^{-1}\dot\phi=R_\theta v.$$
	Consequently, the Lie algebra generating vector fields $g_1,g_2,g_3$ are transformed as follows:
	\begin{equation}\label{vf}
		\begin{aligned}
	g_1&=\partial_{x_1}+\sin(x_4-\tfrac{2\pi}{3})\partial_{x_4}+\sin(x_5)\partial_{x_5}+\sin(x_6+\tfrac{2\pi}{3})\partial_{x_6},\\
	g_2&=\partial_{x_2}-\cos(x_4-\tfrac{2\pi}{3})\partial_{x_4}-\cos(x_5)\partial_{x_5}-\cos(x_6+\tfrac{2\pi}{3})\partial_{x_6},\\
	g_3&=\partial_{x_3}-(1+\cos(x_4))\partial_{x_4}-(1+\cos(x_5))\partial_{x_5}-(1+\cos(x_6))\partial_{x_6}.
	\end{aligned}
	\end{equation}
	We shall use this form for the sake of simplicity.
	Furthermore, to demonstrate the effects of the Lie algebra motions, we calculate the vector fields given by the Lie brackets of $g_1,g_2,g_3$ evaluated at 0 and denote them by $g_{4}=[g_1,g_2],g_{5}=[g_2,g_3]$ and $g_{6}=[g_1,g_3]$. Their coordinates with respect to the system \eqref{eq2} is the following:
	\begin{equation}
	\label{liebrackets}
	\begin{aligned}
	g_4&= \partial x_4+\partial x_5+\partial x_6,\\
	g_5&= \sqrt 3\partial x_5-\sqrt 3\partial x_6,\\
	g_6&= 2\partial x_4-\partial x_5-\partial x_6.
	\end{aligned}
	\end{equation}
	Following \cite{Ishikawa} we demonstrate the motions generated by the Lie brackets. Further details of the Lie bracket exact realizations are described in Section \ref{sect7} and can be found in \cite{Ishikawa}. The following figures show the trajectories of the root centre point, vertices and wheels when a particular Lie bracket motion is realized.
		\begin{figure}[!h]
			\includegraphics[bb=111 277 471 560,clip, width=\textwidth]{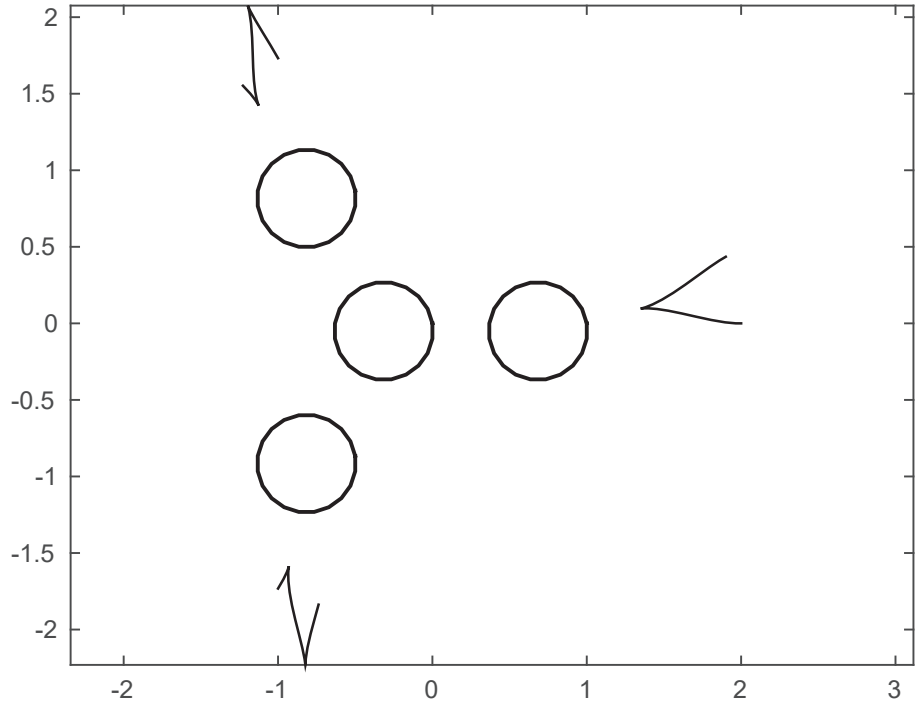}
			\caption{Realization of $g_4$ motion}\label{figg4}
		\end{figure}
	Note that the trajectories on Fig. \ref{figg4} read that the root stays put and the angles represented by the coordinates $x_4,x_5,x_6$ change, which is obvious from approximately equal dislocation of the wheel points at the end of the motion. Considering the vector field $g_4$ at 0 one finds that the angles should change proportionally to 1:1:1.  Similarly, Fig. \ref{figg5} demonstrates the Lie bracket $g_5$ motion and clearly the trajectories represent the effect that the root moves along the $x$--axis and the angles change proportionally to 1:0:-1. Finally, Fig. \ref{figg6} shows $g_6$ realization which reads that the root moves along the $y$--axis and the angles change proportionally to -2:1:1.

	\begin{figure}[!ht]
		\centering
		\begin{minipage}[b]{0.45\textwidth}
			\includegraphics[bb=111 277 471 560,clip, width=\textwidth]{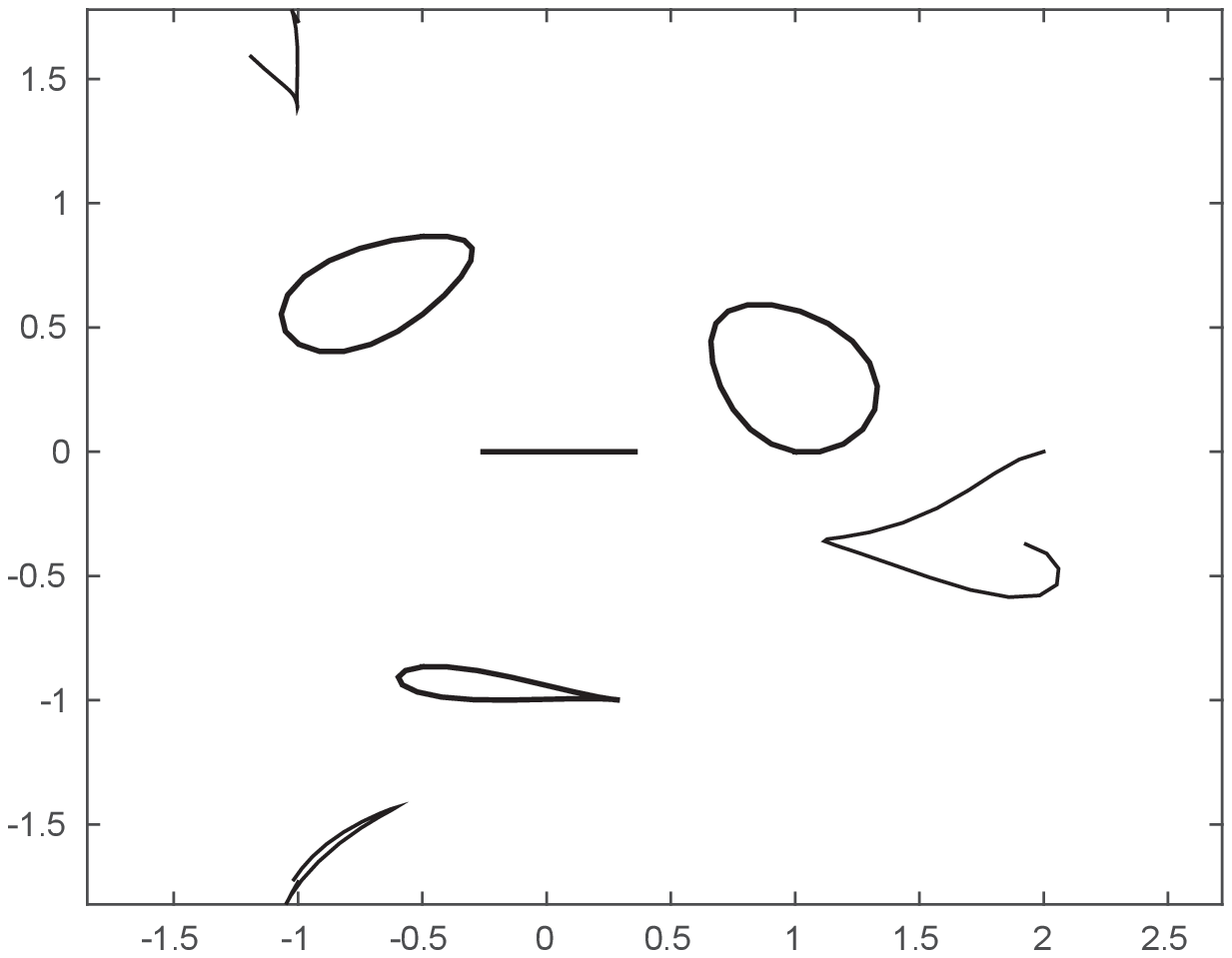}
			\caption{Realization of $g_5$ motion}\label{figg5}
		\end{minipage}
		\hfill
		\begin{minipage}[b]{0.45\textwidth}
			\includegraphics[bb=111 277 471 560,clip,width=\textwidth]{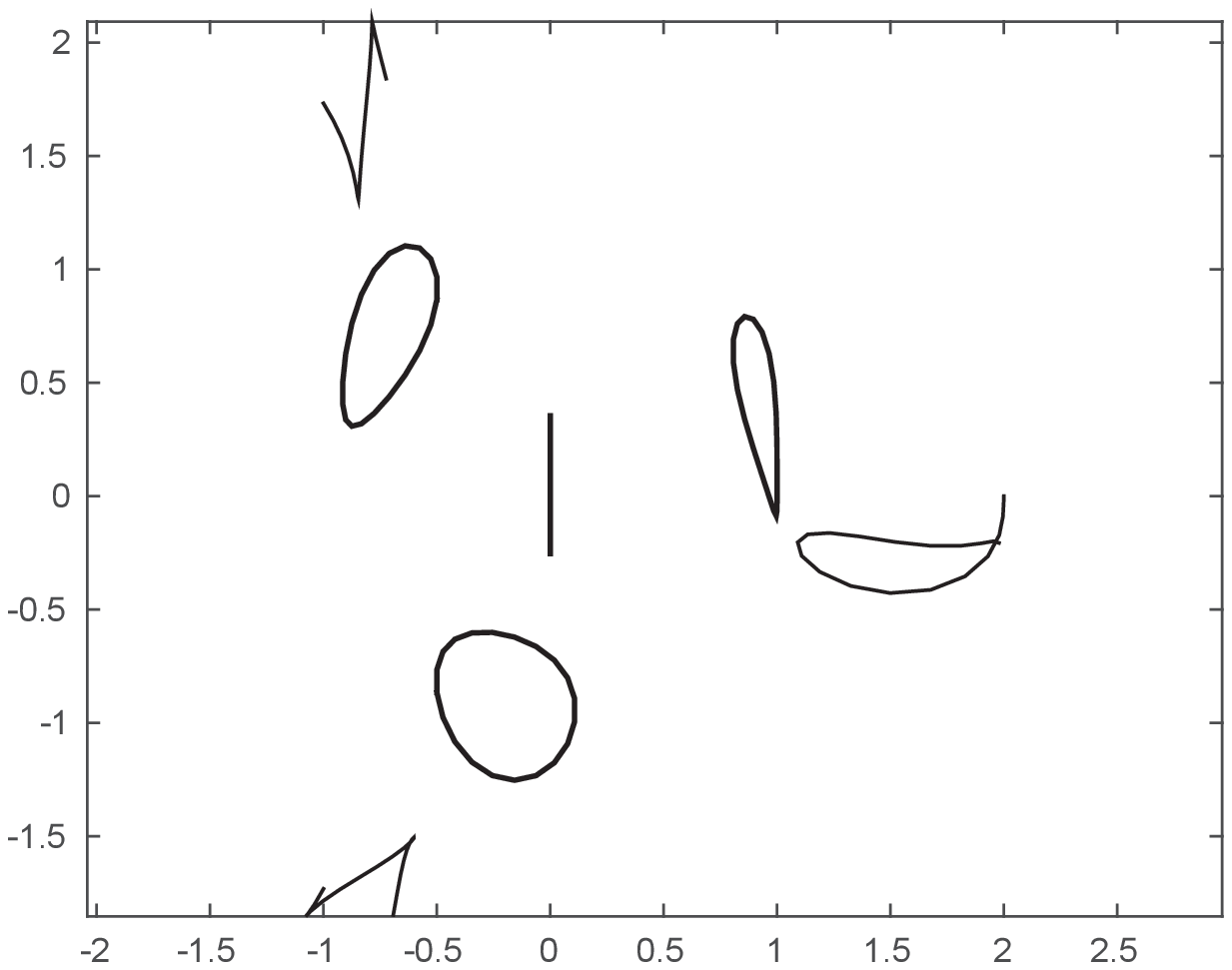}
			\caption{Realization of $g_6$ motion}\label{figg6}
		\end{minipage}
	\end{figure}

	\section{Privileged coordinates}\label{Privcoord}
	A general definition of privileged coordinates is the following, \cite{Jean}, taking into account the notation from Section \ref{Prelim}.
	\begin{definition}\label{privcoorddef}
		A system of privileged coordinates at $p$ is a system of local coordinates $(y_1,...,y_n)$ such that $\mathrm{ord}_p(y_j)=w_j$ for $j=1,...,n.$
	\end{definition}
	In our particular case the configuration space of the trident snake robot is a 6--dimensional manifold $M$ with the coordinate functions denoted by $$(x,y,\theta,\phi_1,\phi_2,\phi_3)=:(x_1,x_2,x_3,x_4,x_5,x_6).$$ Let the basis of a vector space $T_pM$ be denoted by $$(\partial_{x_1},\partial_{x_2},\partial_{x_3}, \partial_{x_4},\partial_{x_5},\partial_{x_6}), p\in M$$ and let us consider three vector fields $g_1,g_2,g_3$ in the form \eqref{vf} which determine a distribution in $TM,$ and we add their Lie brackets $g_4, g_5, g_6,$ see \eqref{liebrackets}. Note that this establishes a filtration of type $(3,6)$ on $TM.$
	The first question is what is the exact form of a coordinate transformation $x:=(x_1,x_2,x_3,x_4,x_5,x_6)\rightsquigarrow (y_1,y_2, y_3,y_4,y_5,y_6)=:y$ such that the condition 
	\begin{equation}\label{cond1}
	\frac{\partial}{\partial y_i}\mid_p=g_i\ \mid_p, \quad i=1,...,6
	\end{equation}
	holds in $p\in M.$ 
	Let us denote by $[g^i_k]_y$ the $i$--th coordinate of a vector $g_k$ in the coordinate system $y$ and by $e^i$ a 6--dimensional vector with coordinates $e^i_j=0$ for $i\neq j$ and $e^i_j=1$ for $i=j,
	i,j\in\{1,...,6\}.$ Then e.g. $[g^1_1]_x=1,\ [g^2_1]_x=0,\ [g_1^3]_x=0,\ [g_1^4]_x=\sin(x_3+x_6)$ etc. and
	the condition \eqref{cond1} reads $[g_i]_y=e^i.$
	Employing the Einstein summation convention, i.e. summing over $j$ ranging from $1$ to $6,$ the transformation law for vector fields under the coordinate change $x\rightsquigarrow y$ reads
	$$\left[g^i_k\right]_y=\frac{\partial y_i}{\partial x_j}\left[g^j_k\right]_x.$$
	Particularly, 
	in the vector form we have
	$$e^i=\left[g_i\right]_y=
	\begin{pmatrix}
	\frac{\partial y_1}{\partial x_j}\left[g^j_i\right]_x\\
	\frac{\partial y_2}{\partial x_j}\left[g^j_i\right]_x\\
	\vdots\\
	\frac{\partial y_6}{\partial x_j}\left[g^j_i\right]_x
	\end{pmatrix}.$$
	Evaluating all functions at an arbitrary point $p$, for sake of simplicity we choose the point $p=(0,0,0,0,0,0)$, 
	we get a system of 36 linear PDEs with respect to $\frac{\partial y_i}{\partial x_j}$ with constant coefficients. We split the system into 6 groups, each containing 6 equations for a particular $y_i,$ determine the inverse matrix and continue by integration. Clearly, at an arbitrary $p\in M$ the desired transformation $x\rightsquigarrow y$ will be linear, in our case it will be given by
	$$\begin{pmatrix}
	y_1\\y_2\\y_3\\y_4\\y_5\\y_6
	\end{pmatrix}=
	\begin{pmatrix}
	\phantom{-}1 & \phantom{-}0 & \phantom{-}0 & \phantom{-}0 & 0 & \phantom{-}0\\
	\phantom{-}0 & \phantom{-}1 & \phantom{-}0 & \phantom{-}0 & 0 & \phantom{-}0\\
	\phantom{-}0 & \phantom{-}0 & \phantom{-}1 & \phantom{-}0 & 0 & \phantom{-}0\\
	\phantom{-}\tfrac{\sqrt 3}{2} & \phantom{-}\tfrac{1}{2} & -2 & 1 & -1 & \phantom{-}\sqrt 3\\
	\phantom{-}0 & -1 & -2 & 1 & \phantom{-}2 & \phantom{-}0\\
	-\tfrac{\sqrt{3}}{2} & \phantom{-}\tfrac{1}{2} & -2 & 1 & -1 & -\sqrt 3
	\end{pmatrix}
	\begin{pmatrix}
	x_1\\x_2\\x_3\\x_4\\x_5\\x_6
	\end{pmatrix}
	$$
	The coordinates $y=(y_1,y_2, y_3,y_4,y_5,y_6)$ are clearly the privileged ones. 
	
	\section{Nilpotent Approximation}\label{NilApprox}
	We proceed according to Bella\"{\i}che\textquoteright s algorithm. Note that in the sequel we use the first two steps only due to the fact that in our filtration (3,6) of $T_pM$ the weights at $p$ are 1 and 2 and thus no further modification of the coordinate system is needed, see \cite{Jean} for a detailed explanation and proof.
	Let us consider the vector fields $g_1,g_2,g_3$ from Section \ref{Privcoord} expressed in the privileged coordinate system $y=(y_1,y_2, y_3,y_4,y_5,y_6)$. 
	
	Vector fields $g_i$ are of order $\geq -1$ and thus generally their Taylor expansion is of the form:
	$$g_i(y)\sim \sum_{\alpha,j}a_{\alpha,j}y^\alpha\partial_{y_j},$$
	where $\alpha=(\alpha_1,...,\alpha_n)$ is a multiindex. Furthermore, if we define a weighted degree of the monomial $y^\alpha=y_1^{\alpha_1}\cdots y_n^{\alpha_n}$ to be $w(\alpha)=w_1\alpha_1+\cdots w_n\alpha_n,$ then $w(\alpha)\geq w_j-1$ if $a_{\alpha,j}\neq 0$. Recall that $w_j=\mathrm{ord}_p(y_j)$ from Definition \ref{privcoorddef} and in our particular case the coordinate weights are $(1,1,1,2,2,2)$. 
	Grouping together the monomial vector fields of the same weighted degree we express $g_i, i=1,2,3$ as a series
	$$g_i=g_i^{(-1)}+g_i^{(0)}+g_i^{(1)}+\cdots,$$
	where $g_i^{(s)}$ is a homogeneous vector field of degree $s$. Note that this means that the $\partial_{y_1},\partial_{y_2}$ and $\partial_{y_3}$ coordinate functions of $g_1^{(-1)}, g_2^{(-1)}$ and $g_3^{(-1)}$ are formed by constants and the $\partial_{y_4},\partial_{y_5}$ and $\partial_{y_6}$ coordinate functions are linear polynomials in $y_1, y_2, y_3.$ Then the following proposition holds, \cite{Jean}:
	
	\begin{proposition}\label{nilp}
		Set $\hat g_i=g_i^{(-1)}, i=1,2,3.$ The family of vector fields $(\hat g_1,\hat g_2,\hat g_3)$ is a first order approximation of $(g_1,g_2,g_3)$ at 0 and generates a nilpotent Lie algebra of step $r=1$, i.e. all brackets of length greater than 1 are zero.
	\end{proposition}
	In our case, we obtain the following vector fields:
	\begin{align*}
	\hat g_1&=  \partial_{y_1}-\frac{y_2}{2}\partial_{y_4}+(-\frac{y_2}{2}-y_3)\partial_{y_5}-\frac{y_1}{2}\partial_{y_6},\\
	\hat g_2&= \partial_{y_2}+\frac{y_1}{2}\partial_{y_4}-\frac{y_1}{2}\partial_{y_5}+(\frac{y_2}{2}-y_3)\partial_{y_6},\\
	\hat g_3 &=\partial_{y_3}.
	\end{align*}
	The family $(\hat g_1,\hat g_2,\hat g_3)$ is the nilpotent approximation of $(g_1, g_2,g_3)$ at 0 associated with the coordinates $y.$ 
	The remaining three vector fields are generated by Lie brackets of $(\hat g_1,\hat g_2,\hat g_3)$ due to the second part of Proposition \ref{nilp}. Note that due to linearity of the three latter coordinates of $(\hat g_1,\hat g_2,\hat g_3)$, the coordinates of $(\hat g_4,\hat g_5,\hat g_6)$ must be constant. We get
	\begin{align*}
	\hat g_4&= \partial_{y_4}\\
	\hat g_5&= \partial_{y_5}\\
	\hat g_6&= \partial_{y_6}.
	\end{align*}
	\section{Lie bracket motion effects}\label{sect7}
\begin{figure}[!h]
		\centering
		\includegraphics[bb=111 277 471 560,scale=.7]{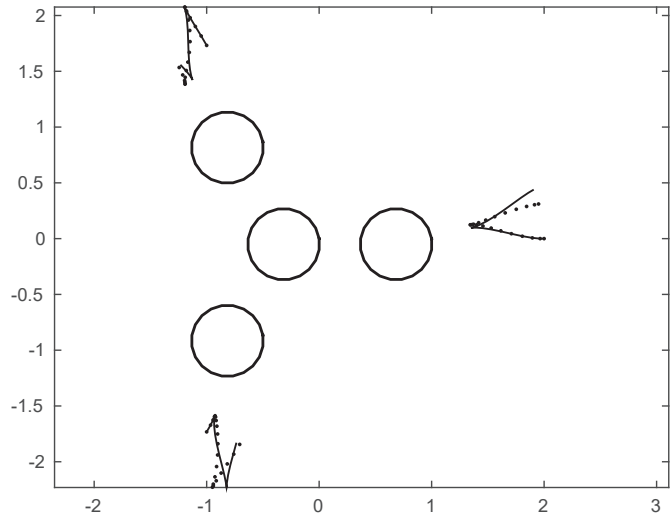}
		\caption{Comparison of $g_4$ motions}\label{fig12}
	\end{figure}

	 In the following, we compare the effect of the Lie bracket motions in the original coordinate system and in the nilpotent approximation. To do so we follow the structure of \cite{Ishikawa}, yet to compare the vector fields in the same coordinate system, the inverse transformation must be applied first and the evaluation of the vector fields effects must be done consequently. Note that the vector fields $(\hat g_1,\hat g_2,\hat g_3,\hat g_4,\hat g_5,\hat g_6)$ in $(x_1,x_2,x_3,x_4,x_5,x_6)$ coordinates are of the form
	 \begin{equation*}
	 \begin{aligned}
	 \hat g_1 &= \partial_{x_1}-(x_2+x_3)\partial_{x_4}-(\tfrac{\sqrt 3x_1}{4}+\tfrac{x_2}{4}-\tfrac{x_3}{2}-\tfrac{\sqrt 3}{2})\partial_{x_5}+\\
	 &\phantom{==}+(\tfrac{\sqrt 3x_1}{2}-
	 \tfrac{x_2}{4}+\tfrac{x_3}{2}-\tfrac{\sqrt 3}{2})\partial_{x_6},\\
	 \hat g_2 &= \partial_{x_2}-\partial_{x_4}+(\tfrac{3x_1}{4}+\tfrac{\sqrt 3x_2}{4}-\tfrac{\sqrt 3x_3}{2}+\tfrac{1}{2})\partial_{x_5}+(\tfrac{3x_1}{4}-\tfrac{\sqrt 3x_2}{4}+\tfrac{\sqrt 3x_3}{2}+\tfrac{1}{2})\partial_{x_6},\\
	 \hat g_3 &= \partial_{x_3}-2\partial_{x_4}-2\partial_{x_5}-2\partial_{x_6}\\
	 \hat g_4 &= \partial_{x_4}+\partial_{x_5}+\partial_{x_6},\\
	 \hat g_5 &= -\tfrac{\sqrt 3}{2}\partial_{x_5}+\tfrac{\sqrt 3}{2}\partial_{x_6},\\
	 \hat g_6 &= -\partial_{x_3}+\tfrac{1}{2}\partial_{x_5}+\tfrac{1}{2}\partial_{x_6}.
	 \end{aligned}
	 \end{equation*}
	Note that the Lie bracket motions at 0 correspond exactly to the original ones. Anyway, to perform the Lie bracket motions we apply a periodic input, i.e. for the vector fields $\hat g_4=[\hat g_1,\hat g_2], \hat g_5=[\hat g_1,\hat g_3], \hat g_6=[\hat g_2,\hat g_3],$ respectively, the input  
	\begin{align}
	v_1(t)&=(-A\omega\sin\omega t, A\omega\cos\omega t,0)\label{periodic_input4}\\
	v_2(t)&=(0,-A\omega\sin\omega t, A\omega\cos\omega t)\label{periodic_input6}\\
	v_3(t)&=(-A\omega\sin\omega t, 0,A\omega\cos\omega t)\label{periodic_input5}
	\end{align}
	is applied, because, according to \cite{Ishikawa}, the Lie bracket of a pair of vector fields corresponds to the direction of a displacement in the state space as a result of a periodic input with sufficiently small amplitude $A$, i.e. the bracket motions are generated by periodic combination of the vector controlling fields. In Fig. \ref{fig12}, there is a comparison of the $g_4$ motion realized by the periodic input in $x_1,...,x_6$ coordinates (dotted line) and in nilpotent approximation. 
	
\begin{figure}[!ht]
		\centering
		\begin{minipage}[b]{0.45\textwidth}
			\includegraphics[bb=111 277 471 560,clip, width=\textwidth]{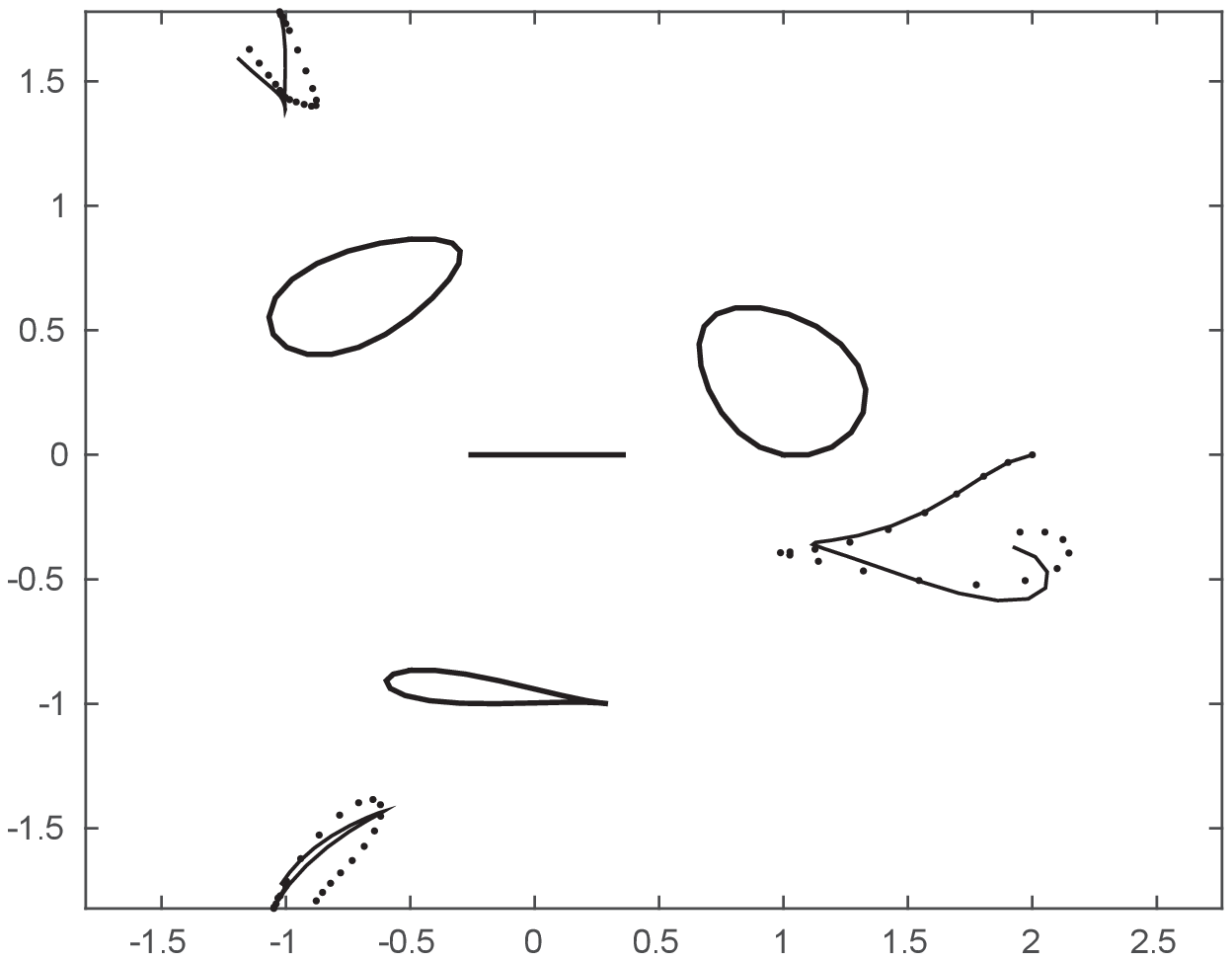}
			\caption{Comparison of $g_5$ motions}\label{fig13}
		\end{minipage}
		\hfill
		\begin{minipage}[b]{0.45\textwidth}
			\includegraphics[bb=111 277 471 560,clip,width=\textwidth]{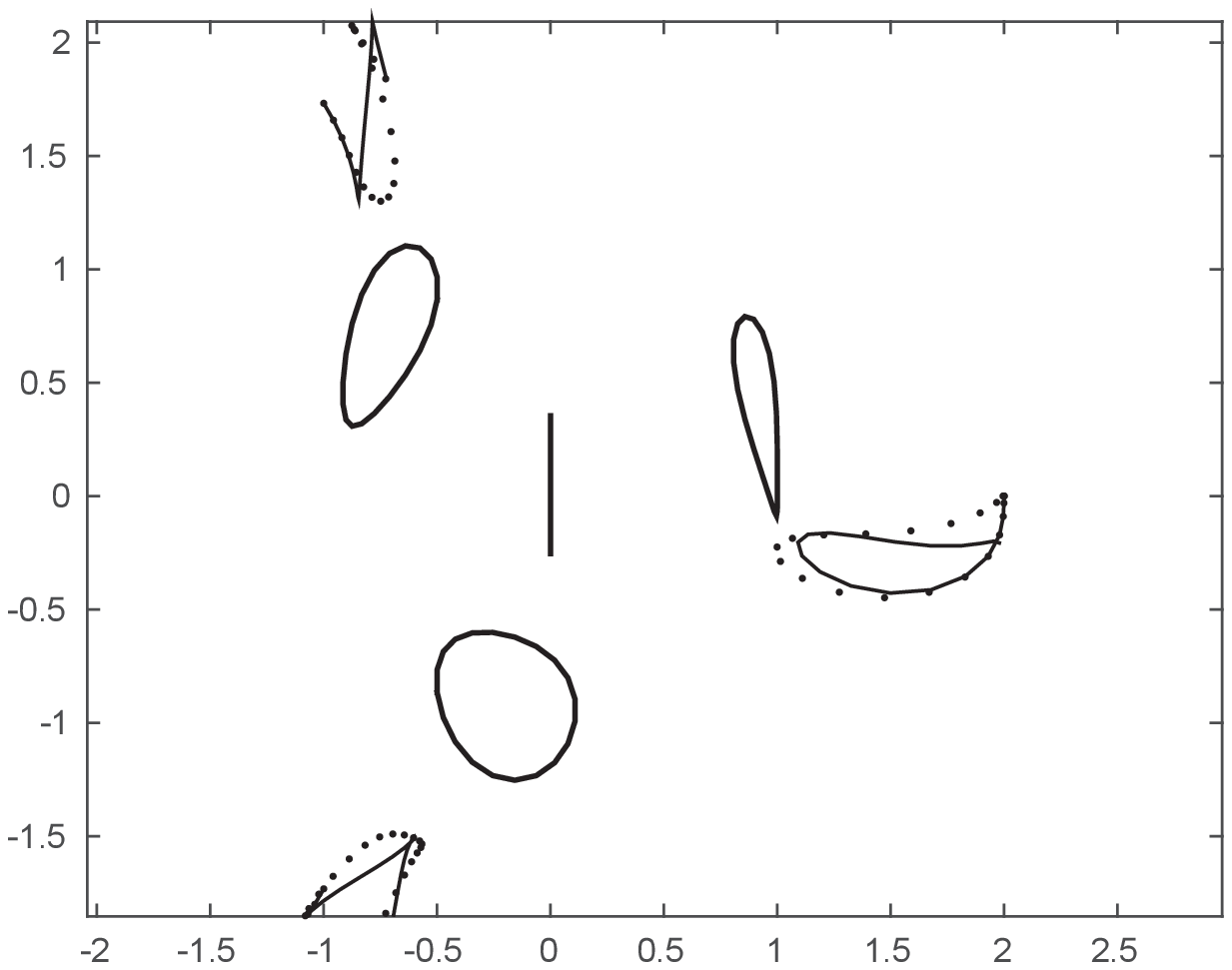}
			\caption{Comparison of $g_6$ motions}\label{fig23}
		\end{minipage}
	\end{figure}

	Fig. \ref{fig13} and \ref{fig23} show the comparison of $g_5$ and $g_6$ motions, respectively. Note that the lines represent the trajectories of the appropriate wheel and thus the accuracy of the motion in real space is pictured.

	\section{Conclusions}
	We presented a calculation of a nilpotent approximation of the family of vector fields corresponding to the controlling distribution of a trident snake robot. Such an approximation is valuable not only for the calculational complexity reasons but also from the theoretical point of view as the nilpotency  simplifies the model for further theoretical considerations significantly. We showed that even from the practical point of view this approximation is good as the deviation from the exact model control is minimal. More precisely, we checked that at 0 the Lie brackets of the original controlling vector fields and of the approximated ones coincide and, furthermore, if their realization by the periodic input is considered, the deviations depicted in Figures \ref{fig12}, \ref{fig13}, \ref{fig23} are minimal. Finally let us claim that the error in control leads to the violation of the non--holonomic conditions and thus the wheels slip a bit, yet the benefits of the nilpotent approximation prevail.

	
\email{hrdina@fme.vutbr.cz,navrat.a@fme.vutbr.cz, vasik@fme.vutbr.cz}
\end{document}